% ...
\documentclass[a4paper,12pt]{article}

\usepackage{graphicx}
\usepackage{fullpage}

\usepackage{amsthm, amsmath, amsfonts, amssymb}

\theoremstyle{plain}
\newtheorem{definition}{Definition}
\theoremstyle{plain}
\newtheorem{theorem}[definition]{Theorem}
\newtheorem{lemma}[definition]{Lemma}
\newtheorem{corollary}[definition]{Corollary}
\theoremstyle{definition}

\def\copyrightheading{~}
\newcommand{\fcstitle}[1]{{\LARGE {#1}}}
\def\authorfont{\large}
\def\addressfont{\normalfont}
\def\smalllineskip{\smallskip}
\def\textlineskip{\medskip}

\def\epsilon{\varepsilon}

\def\dtv{d_{\textup{TV}}}

\def\Prob{\textnormal{Pr}}

\def\Pk{P^{[k]}}

\def\P1{P^{[1]}}

\def\Pgrid{P_{\textup{grid}}}

\def\M{\mathcal{M}}
\def\gridscan{\mathcal{M}_{\textup{grid}}}
\def\scan{\mathcal{M}_{\rightarrow}}

\def\Mix{\textup{Mix}}

\def\pone{0.283}
\def\ptwo{0.079}
\def\pthree{0.051}
\def\pfour{0.079}
\def\alphasum{0.984}

\def\ponesmall{0.379}
\def\ptwosmall{0.107}
\def\pthreesmall{0.050}
\def\pfoursmall{0.107}
\def\alphasumsmall{1.286}

\def\ponemedium{0.3671}
\def\pthreemedium{0.0298}
\def\pfourmedium{0.0997}
\def\psixmedium{0.0174}
\def\alphasummedium{1.028}

\def\ponecornerbig{0.3537}
\def\pthreecornerbig{0.0245}
\def\psevencornerbig{0.0245}
\def\pninecornerbig{0.0071}
\def\ponemiddlebig{0.0838}
\def\pthreemiddlebig{0.0838}
\def\psevenmiddlebig{0.0138}
\def\pninemiddlebig{0.0138}
\def\alphasumbig{1.0148}

\def\scanURL{\texttt{http://www.csc.liv.ac.uk/\nolinebreak[3]$\sim$markus/\nolinebreak[3]systematicscan/}}

\begin{document}
\copyrightheading



\begin{center}

\fcstitle{A Systematic Scan for $7$-colourings of the Grid}

\vspace{24pt}

{\authorfont Markus Jalsenius}

\vspace{2pt}

\smalllineskip
{\addressfont Department of Computer Science, University of Liverpool\\
Ashton Street, Liverpool, L69 3BX, United Kingdom}

\vspace{10pt}

and

\vspace{10pt}

{\authorfont Kasper Pedersen}

\vspace{2pt}
\smalllineskip
{\addressfont Department of Computer Science, University of Liverpool\\
Ashton Street, Liverpool, L69 3BX, United Kingdom}

\vspace{20pt}

\end{center}


\begin{abstract}
We study the mixing time of a systematic scan Markov chain for sampling
from the uniform distribution on proper $7$-colourings of a finite
rectangular sub-grid of the infinite square lattice, the grid. A systematic
scan Markov chain cycles through finite-size subsets of vertices in
a deterministic order and updates the colours assigned to the vertices
of each subset. The systematic scan Markov chain that we present cycles
through subsets consisting of 2$\times$2 sub-grids and updates the
colours assigned to the vertices using a procedure known as heat-bath.
We give a computer-assisted proof that this systematic scan Markov
chain mixes in $O(\log n)$ scans, where $n$ is the size of the rectangular
sub-grid. We make use of a heuristic to compute required couplings
of colourings of 2$\times$2 sub-grids. This is the first time the
mixing time of a systematic scan Markov chain on the grid has been
shown to mix for less than $8$ colours. We also give partial results
that underline the challenges of proving rapid mixing of a systematic
scan Markov chain for sampling $6$-colourings of the grid by considering
2$\times$3 and 3$\times$3 sub-grids.
\end{abstract}

\textlineskip

\section{Introduction}

\pagenumbering{arabic}This paper is concerned with sampling from
the uniform distribution, $\pi$, on the set of proper $q$-colourings
of a finite-size rectangular grid. A \emph{$q$-colouring} of a graph
is an assignment of a colour from a finite set of $q$ distinct colours
to each vertex and we say that a colouring is a \emph{proper} \emph{colouring}
if no two adjacent vertices are assigned the same colour. Proper $q$-colourings
of the grid correspond to the zero-temperature anti-ferromagnetic
$q$-state Potts model on the square lattice, a model of significant
importance in statistical physics (see for example Salas and Sokal~\cite{salas-sokal}). 

Sampling from $\pi$ is computationally challenging, however it remains
an important task and it is frequently carried out in experimental
work by physicists by simulating some suitable random dynamics that
converges to $\pi$. Ensuring that a dynamics converges to $\pi$
is generally straight forward, but obtaining good upper bounds on
the number of steps required for the dynamics to become sufficiently
close to $\pi$ is a much more difficult problem. Physicists are at
times forced to ``guess'' (using some heuristic methods) the number
of steps required for their dynamics to be sufficiently close to the
uniform distribution in order to carry out their experiments. By establishing
rigorous bounds on the convergence rates (\emph{mixing time}) of these
dynamics computer scientists can provide underpinnings for this type
of experimental work and also allow a more structured approach to
be taken. 

Providing bounds on the mixing time of Markov chains is a well-studied
problem in theoretical computer science. However, the types of Markov chains
frequently considered by computer scientists do not always
correspond to the dynamics usually used in the experimental work by
physicists. In computer science, the mixing time of various types
of \emph{random update Markov chains} have been frequently analysed;
notably on the grid by Achlioptas, Molloy, Moore and van Bussel~\cite{grid_ach}
and Goldberg, Martin and Paterson~\cite{ssm}. We say that a Markov
chain on the set of colourings is a random update Markov chain when
one step of the the process consists of randomly selecting a set of
vertices (often a single vertex) and updating the colours assigned
to those vertices according to some well-defined distribution induced
by $\pi$. Experimental work is, however, often carried out by cycling
through and updating the vertices (or subsets of vertices) in a deterministic
order. This type of dynamics has recently been studied by computer
scientists in the form of \emph{systematic scan Markov chains} (\emph{systematic
scan} for short). For results regarding systematic scan see for instance
Dyer, Goldberg and Jerrum~\cite{systematic_scan,dobrushin_scan}
and Pedersen~\cite{dobrushin_paper} although these papers are not
considering the grid specifically. It is important to note that systematic
scan remains a random process since the method used to update the
colour assigned to the selected set of vertices is a randomised procedure
drawing from some well-defined distribution induced by $\pi$. 

In Section~\ref{sec:coupling} we present a computer assisted proof
that systematic scan mixes rapidly when considering $7$-colourings
of the grid. Previously eight was the least number of colours for
which systematic scan on the grid was known to be rapidly mixing,
due to Pedersen~\cite{dobrushin_paper}, a result which we hence
improve on in this paper. We will make use of a recent result by Pedersen~\cite{dobrushin_paper}
to prove rapid mixing of systematic scan by bounding the influence
\emph{on} a vertex (note that the literature traditionally talks about
sites rather than vertices). We will provide bounds on this influence
parameter by using a heuristic to mechanically construct sufficiently
good \emph{couplings} of proper colourings of a 2$\times$2~sub-grid.
We will hence use a heuristic based computation in order to establish
a rigorous result about the mixing time of a systematic scan Markov
chain. Finally, in Section~\ref{sec:Partial-results}, we consider
the possibility of proving rapid mixing of systematic scan for $6$-colourings
of the grid by increasing the size of the sub-grids. We give lower
bounds on the appropriate influence parameter that imply that the
proof technique we employ does not imply rapid mixing of systematic
scan for $6$-colourings of the grid when using 2$\times$2, 2$\times$3
and~3$\times$3 sub-grids.

\subsection{Preliminaries and statement of results}

Let $Q=\{1,\dots,7\}$ be the set of colours and $V=\{1,\ldots,n\}$
the set of vertices of a finite rectangular grid $G$ with toroidal
boundary conditions. Working on the torus is common practice as it
avoids treating several technicalities regarding the vertices on the
boundary of a finite grid as special cases and hence lets us present
the proof in a more ``clean'' way. We point out however that these
technicalities are straightforward to deal with (more on this in Section~\ref{sec:mixing}).
We formally say that
a colouring $\sigma$ of $G$ is a function from $V$ to $Q$. Let
$\Omega^{+}$ be the set of all colourings of $G$ and $\Omega$ be
the set of all proper $q$-colourings. Then the distribution $\pi$,
described earlier, is the uniform distribution on $\Omega$. If $\sigma\in\Omega^{+}$
is a colouring and $j\in V$ is a vertex then $\sigma_{j}$ denotes
the colour assigned to vertex $j$ in colouring $\sigma$. Furthermore,
for a subset of vertices $\Lambda\subseteq V$ and a colouring $\sigma\in\Omega^{+}$
we let $\sigma_{\Lambda}$ denote the colouring of the vertices in
$\Lambda$ under $\sigma$. For each vertex $j\in V$, let $S_{j}$
denote the set of pairs $(\sigma,\tau)\in\Omega^{+}\times\Omega^{+}$
of colourings that only differ on the colour assigned to vertex $j$,
that is $\sigma_{i}=\tau_{i}$ for all $i\neq j$. 

Let $\M$ be a Markov chain with state space $\Omega^{+}$ and stationary
distribution $\pi$. Suppose that the transition matrix of $\M$ is
$P$. Then the mixing time from an initial colouring $\sigma\in\Omega^{+}$
is the number of steps, that is applications of $P$, required for
$\M$ to become sufficiently close to $\pi$. Formally the mixing
time of $\M$ from an initial colouring $\sigma\in\Omega^{+}$ is
defined, as a function of the deviation $\epsilon$ from stationarity,
by
\begin{equation}
\Mix_{\sigma}(\M,\epsilon)=\min\{ t>0\;:\;\dtv(P^{t}(\sigma,\cdot),\pi)\leq\epsilon\},
\end{equation}
where
\begin{equation}
\dtv(\theta_{1},\theta_{2})=\frac{1}{2}\sum_{i}|\theta_{1}(i)-\theta_{2}(i)|=\max_{A\subseteq\Omega^{+}}|\theta_{1}(A)-\theta_{2}(A)|
\end{equation}
is the total variation distance between two distributions $\theta_{1}$
and $\theta_{2}$ on $\Omega^{+}$. The mixing time $\Mix(M,\epsilon)$
of $\M$ is then obtained my maximising over all possible initial
colourings
\begin{equation}
\Mix(\M,\epsilon)=\max_{\sigma\in\Omega^{+}}\Mix_{\sigma}(\M,\epsilon).
\end{equation}
We say that $\M$ is \emph{rapidly mixing} if the mixing time of $\M$
is polynomial in $n$ and $\log(\epsilon^{-1})$.

We will make use of a recent result by Pedersen~\cite{dobrushin_paper}
to study the mixing time of a systematic scan Markov chain for $7$-colourings
of the grid using \emph{block updates}. We need the following notation
in order to define our systematic scan Markov chain. Define the following
set $\Theta=\{\Theta_{1},\dots,\Theta_{m}\}$ of $m$ blocks. Each
block $\Theta_{k}\subseteq V$ is a 2$\times$2~sub-grid and $m$
is the smallest integer such that $\bigcup_{k=1}^{m}\Theta_{k}=V$.
For any block $\Theta_{k}$ and a pair of colourings $\sigma,\tau\in\Omega^{+}$
we write {}``$\sigma=\tau$ on $\Theta_{k}$'' if $\sigma_{i}=\tau_{i}$
for each $i\in\Theta_{k}$ and similarly {}``$\sigma=\tau$ off $\Theta_{k}$''
if $\sigma_{i}=\tau_{i}$ for each $i\in V\setminus\Theta_{k}$. We
also let $\partial\Theta_{k}$ denote the set of vertices in $V\setminus\Theta_{k}$
that are adjacent to some vertex in $\Theta_{k}$, and we will refer
to $\partial\Theta_{k}$ as the \emph{boundary} of $\Theta_{k}$.
Note from our previous definitions that $\sigma_{\partial\Theta_{k}}$
denotes the colouring of the boundary of $\Theta_{k}$ under a colouring
$\sigma\in\Omega^{+}$. We will refer to $\sigma_{\partial\Theta_{k}}$
as a \emph{boundary colouring}. Finally we say that a $7$-colouring
of the 2$\times$2~sub-grid $\Theta_{k}$ \emph{agrees} with a boundary
colouring $\sigma_{\partial\Theta_{k}}$ if (1) no adjacent sites
in $\Theta_{k}$ are assigned the same colour and (2) each vertex
$j\in\Theta_{k}$ is assigned a colour that is different to the colours
of all boundary vertices adjacent to $j$. 

For each block $\Theta_{k}$ and colouring $\sigma\in\Omega^{+}$
let $\Omega_{k}(\sigma)$ be the subset of $\Omega^{+}$ such that
if $\sigma'\in\Omega_{k}(\sigma)$ then $\sigma'=\sigma$ off $\Theta_{k}$
and $\sigma'_{\Theta_{k}}$ agrees with $\sigma_{\partial\Theta_{k}}$.
Let $\pi_{k}(\sigma)$ be the uniform distribution on $\Omega_{k}(\sigma)$.
We then define $\Pk$ to be the transition matrix on the state space
$\Omega^{+}$ for performing a so-called \emph{heat-bath} move on
$\Theta_{k}$. A heat-bath move on a block $\Theta_{k}$, given a
colouring $\sigma\in\Omega^{+}$, is performed by drawing a new colouring
from the distribution $\pi_{k}(\sigma)$. Note in particular that
applying $\Pk$ to a colouring $\sigma\in\Omega^{+}$ results in a
colouring $\sigma^{\prime}\in\Omega^{+}$ such that $\sigma'=\sigma$
off $\Theta_{k}$ and the colouring $\sigma_{\Theta_{k}}'$ of $\Theta_{k}$
is proper and agrees with the colouring $\sigma_{\partial\Theta_{k}}'$
of the boundary of $\Theta_{k}$ (which is identical to $\sigma_{\partial\Theta_{k}}$).
We formally define the following systematic scan Markov chain for
$7$-colourings of $G$, which systematically performs heat-bath moves
on 2$\times$2 sub-grids, as follows. It is worth pointing out that
this holds for \emph{any} \emph{ordering} of the set of blocks.

\begin{definition}
\label{def:scan}The systematic scan dynamics for $7$-colourings
of $G$ is a Markov chain $\gridscan$ with state space $\Omega^{+}$
and transition matrix $\Pgrid=\Pi_{k=1}^{m}\Pk$. 
\end{definition}
It can be shown that the stationary distribution of $\gridscan$ is
$\pi$ by considering the construction of $\Pgrid$. It is customary
to refer to one application of $\Pgrid$ (that is updating each block
once) as one \emph{scan}. One scan takes $\sum_{k}|\Theta_{k}|$ \emph{vertex
updates} and by construction of $\Theta$ this sum is clearly of order
$O(n)$. 

We will prove the following theorem and point out that this is the
first proof of rapid mixing of systematic scan for $7$-colourings
on the grid. 

\begin{theorem}
\label{thm:gridmix}Let $\gridscan$ be the Markov chain from Definition~\ref{def:scan}
on $7$-colourings of $G$. Then the mixing time of $\gridscan$ is
\begin{equation}
\Mix(\gridscan,\epsilon)\leq63\log(n\epsilon^{-1}).
\end{equation}

\end{theorem}

\subsection{Context and related work}

We now provide an overview of previous achievements for colourings
of the grid. Previously it was known that systematic scan for $q$-colourings
on general graphs with maximum vertex degree $\Delta$ mixes in $O(\log n)$
scans when $q\geq2\Delta$ due to Pedersen~\cite{dobrushin_paper}.
That result is a hand-proof and uses block updates that updates the
colour at each endpoint of an edge during each step. Earlier Dyer
et al.~\cite{dobrushin_scan} had shown that a single-site systematic
scan Markov chain (where one vertex is updated at a time) mixes in
$O(\log n)$ scans when $q>2\Delta$ and in $O(n^{2}\log n)$ scans
when $q=2\Delta$. It is hence well-established that systematic scan
is rapidly mixing for $q$-colourings of the grid when $q\geq8$ but
nothing has been known about the mixing time for smaller $q$. The
results of both Pedersen~\cite{dobrushin_paper} and Dyer et al.~\cite{dobrushin_scan}
bound the mixing time by studying the influence \emph{on} a vertex.
We will use that technique in this paper as well, however we will
construct the required couplings using a heuristic. We defer the required
definitions to Section~\ref{sec:mixing} which also contains the
proof of Theorem~\ref{thm:gridmix}.

Recent results have revealed that, in a single-site setting, one is
not restricted use the total influence on a vertex when analysing the
mixing time of systematic scan by bounding influence parameters. In a
single-site setting one can define an $n$$\times$$n$-matrix whose
entries are the influences that all vertices have on each other. Hayes~\cite{hayes}
has shown that providing a sufficiently small upper bound on the
spectral gap of this matrix implies rapid mixing of both systematic
scan and random update. Dyer, Goldberg and Jerrum~\cite{matrixnorm} furthermore
showed that an upper bound on any matrix norm also implies rapid
mixing of both types of Markov chains. These techniques are however
not known to apply to Markov chains using block moves. See the PhD
thesis by Pedersen~\cite{thesis} for more comprehensive review of the above
results and for the difficulties in extending them to cover block
dynamics.

As random update Markov chains have received more attention than systematic
scan we also summarise some mixing results of interest regarding $q$-colourings
of the grid (recall that a random update Markov chain selects randomly
a subset of sites to be updated at each step). Achlioptas et al.~\cite{grid_ach}
give a computer-assisted proof of mixing in $O(n\log n)$ updates
when $q=6$ by considering blocks consisting of 2$\times$3~sub-grids.
Our computations are similar in nature to the ones of Achlioptas et
al. however their computations are not sufficient to imply mixing
of systematic scan as we will discuss in due course. More recently
Goldberg, Martin and Paterson~\cite{ssm} gave a hand-proof of mixing
in $O(n\log n)$ updates when $q\geq7$ using the technique of strong
spatial mixing. Previously Salas and Sokal~\cite{salas-sokal} gave
a computer-assisted proof of the $q=7$ case, a result which was also
implied by another computer-assisted result due to Bubley, Dyer and
Greenhill~\cite{bubley_computer-assist} that applies to $4$-regular
triangle-free graphs. Finally it is worth pointing out that, in the
special case when $q=3$, two complementary results of Luby, Randall
and Sinclair~\cite{luby-planar} and Goldberg, Martin and Paterson~\cite{3coloursz2}
give rapid mixing of random update.

\section{Bounding the mixing time of systematic scan\label{sec:mixing}}

This section will contain a proof of Theorem~\ref{thm:gridmix} although
the proof of a crucial lemma, which requires computer-assistance,
is deferred to Section~\ref{sec:coupling}. We will bound the mixing
time of $\gridscan$ by bounding the influence on a vertex, a parameter
which we denote by $\alpha$ and will define formally in due course.
If $\alpha$ is sufficiently small then Theorem~2 from Pedersen~\cite{dobrushin_paper}
implies that any systematic scan Markov chain, whose transition matrices
for updating each block satisfy two simple properties, mixes in $O(\log n)$
scans. For completeness we restate this theorem (Theorem~\ref{thm:main-d} below) and in the statement
we let $\scan$ denote a systematic scan Markov chain whose transition
matrices for each block update satisfy the required properties. 

\begin{theorem}
\label{thm:main-d}
If
$\alpha<1$ then the mixing time of $\scan$ is
\begin{equation}
\Mix(\scan,\epsilon)\leq\frac{\log(n\epsilon^{-1})}{1-\alpha}.
\end{equation}

\end{theorem}
For each block $\Theta_{k}$ the transition matrix $\Pk$ needs to
satisfy the following two properties in order for Theorem~\ref{thm:main-d}
to apply.

\begin{enumerate}
\item If $\Pk(\sigma,\tau)>0$ then $\sigma=\tau$ off $\Theta_{k}$, and
\item $\pi$ is invariant with respect to $\Pk.$
\end{enumerate}
It is pointed out in Pedersen \cite{dobrushin_paper} that if $\Pk$
is a transition matrix performing a heat-bath move then both of these
properties are easily satisfied. Furthermore, it is pointed out that
when $\Omega$ is the set of proper colourings of a graph, then $\pi$
is the uniform distribution on $\Omega$ as we require. Since the
transition matrices $\Pk$ used in the definition of $\gridscan$
perform heat-bath updates we are hence able to use Theorem~\ref{thm:main-d}
to bound the mixing time of $\gridscan$.

We are now ready to formally define the parameter $\alpha$ denoting
the influence on a vertex. For any pair of colourings $(\sigma,\tau)\in S_{i}$
let $\Psi_{k}(\sigma,\tau)$ be a coupling of the distributions induced
by $\Pk(\sigma,\cdot)$ and $\Pk(\tau,\cdot)$, namely $\pi_{k}(\sigma)$
and $\pi_{k}(\tau)$ respectively. We remind the reader that a coupling
of two distributions $\pi_{1}$ and $\pi_{2}$ on state space $\Omega^{+}$
is a joint distribution $\Omega^{+}\times\Omega^{+}$ such that the
marginal distributions are $\pi_{1}$ and $\pi_{2}$. For ease of
reference we also let $p_{j}(\Psi_{k}(\sigma,\tau))$ denote the probability
that a vertex $j\in\Theta_{k}$ is assigned a different colour in
a pair of colourings drawn from some coupling $\Psi_{k}(\sigma,\tau)$.
We then let
\begin{equation}
\rho_{i,j}^{k}=\max_{(\sigma,\tau)\in S_{i}}p_{j}(\Psi_{k}(\sigma,\tau))
\end{equation}
be the influence of $i$ on $j$ under $\Theta_{k}$. Finally the
parameter $\alpha$ denoting the influence on any vertex is defined
as
\begin{equation}
\alpha=\max_{k}\max_{j\in\Theta_{k}}\sum_{i}\rho_{i,j}^{k}.
\end{equation}

Pedersen~\cite{dobrushin_paper} actually defines $\alpha$ with
a weight associated with each vertex, however as we will not use weights
in our proof we have omitted them from the above account. 
So, in order to upper bound $\alpha$ we are required to upper bound
the probability of a discrepancy at each vertex $j\in\Theta_{k}$
under a coupling $\Psi_{k}(\sigma,\tau)$ of the distributions $\pi_{k}(\sigma)$
and $\pi_{k}(\tau)$ for any pair of colourings $(\sigma,\tau)\in S_{i}$
that only differ at the colour of vertex $i$. Our main task is hence
to specify a coupling $\Psi_{k}(\sigma,\tau)$ of $\pi_{k}(\sigma)$
and $\pi_{k}(\tau)$ for each pair of colourings $(\sigma,\tau)\in S_{i}$
and upper bound the probability of assigning a different colour to
each vertex in a pair of colourings drawn from that coupling. 

Consider any block $\Theta_{k}$ and any pair of colourings $(\sigma,\tau)\in S_{i}$
that differ only on the colour assigned to some vertex $i$. Clearly
the distribution on colourings of $\Theta_{k}$, induced by $\pi_{k}(\sigma)$
only depends on the boundary colouring $\sigma_{\partial\Theta_{k}}$.
Similarly, the distribution on colourings of $\Theta_{k}$, induced
by $\pi_{k}(\tau)$ depends only on $\tau_{\partial\Theta_{k}}$.
If $i\not\in\partial\Theta_{k}$ then the distributions on the colourings
of $\Theta_{k}$, induced by $\pi_{k}(\sigma)$ and $\pi_{k}(\tau)$,
respectively, are the same and we let $\Psi_{k}(\sigma,\tau)$ be
the coupling in which any pair of colourings drawn from $\Psi_{k}(\sigma,\tau)$
agree on $\Theta_{k}$. That is, if the pair $(\sigma',\tau')$ of
colourings are drawn from $\Psi_{k}(\sigma,\tau)$ then $\sigma'=\sigma$
off $\Theta_{k}$, $\tau'=\tau$ off $\Theta_{k}$ and $\sigma'=\tau'$
on $\Theta_{k}$. This gives $\rho_{i,j}^{k}=0$ for any $i\not\in\partial\Theta_{k}$
and $j\in\Theta_{k}$. 

We now need to construct $\Psi_{k}(\sigma,\tau)$ for the case when
$i\in\partial\Theta_{k}$. For each $j\in\Theta_{k}$ we need $p_{j}(\Psi_{k}(\sigma,\tau))$
to be sufficiently small in order to avoid $\rho_{i,j}^{k}$ being
too big. If the $\rho_{i,j}^{k}$-values are too big the parameter
$\alpha$ will be too big (that is greater than one) and we cannot
make use of Theorem~\ref{thm:main-d} to show rapid mixing. Constructing
$\Psi_{k}(\sigma,\tau)$ by hand such that $p_{j}(\Psi_{k}(\sigma,\tau))$
is sufficiently small is a difficult task. It is, however, straight
forward to mechanically determine which colourings have positive measure
in the distributions $\pi_{k}(\sigma)$ and $\pi_{k}(\tau)$ for a
given pair of boundary colourings $\sigma_{\partial\Theta_{k}}$ and
$\tau_{\partial\Theta_{k}}$. From these distributions we can then
use some suitable heuristic to construct a coupling that is good enough
for our purposes. We hence need to construct a specific coupling for
each individual pair of colourings differing only at a single vertex.
In order to do this we will make use of the following lemma, which
is proved in Section~\ref{sec:coupling}. 

\begin{figure}
\centering
\includegraphics[scale=0.8]{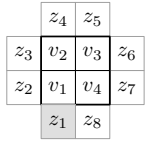}
\caption{\label{fig:labeling-2x2-block}General labeling of the vertices in
a 2$\times$2-block $\Theta_{k}$ and the vertices $\partial\Theta_{k}$
on the boundary of the block.}
\end{figure}

\begin{lemma}
\label{lem:main-lemma}Let $v_{1},\dots,v_{4}$ be the four vertices
in a 2$\times$2-block and $z_{1},\dots,z_{8}$ be the boundary vertices
of the block and let the labeling be as in Figure~\ref{fig:labeling-2x2-block}.
Let $Z$ and $Z'$ be any two $7$-colourings of the boundary vertices
such that $Z$ and $Z'$ agree on each vertex except on $z_{1}$.
Let $\pi_{Z}$ and $\pi_{Z'}$ be the uniform distributions on proper
$7$-colourings of the block that agree with $Z$ and $Z'$, respectively.
For $i=1,\dots,4$ let $p_{v_{i}}(\Psi)$ denote the probability
that the colour of vertex $v_{i}$ differ in a pair of colourings
drawn from a coupling $\Psi$ of $\pi_{Z}$ and $\pi_{Z'}$. Then
there exists a coupling $\Psi$ such that $p_{v_{1}}(\Psi)<\pone$,
$p_{v_{2}}(\Psi)<\ptwo$, $p_{v_{3}}(\Psi)<\pthree$ and $p_{v_{4}}(\Psi)<\pfour$.
\end{lemma}

\begin{figure}
\centering
(a)\includegraphics[scale=0.8]{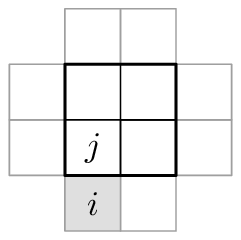}
\hspace{2mm}
(b)\includegraphics[scale=0.8]{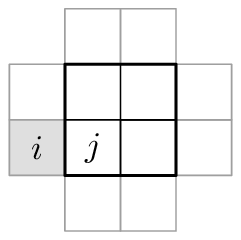}
\hspace{2mm}
(c)\includegraphics[scale=0.8]{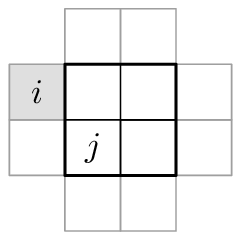}
\hspace{2mm}
(d)\includegraphics[scale=0.8]{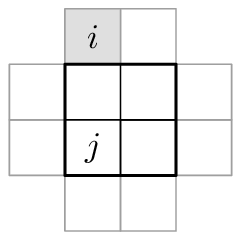}
\medskip

(e)\includegraphics[scale=0.8]{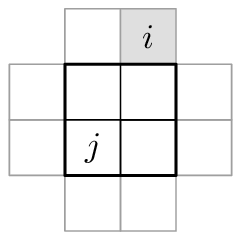}
\hspace{2mm}
(f)\includegraphics[scale=0.8]{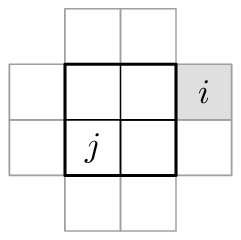}
\hspace{2mm}
(g)\includegraphics[scale=0.8]{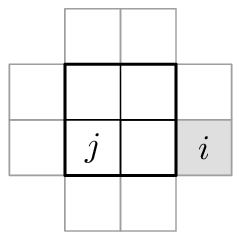}
\hspace{2mm}
(h)\includegraphics[scale=0.8]{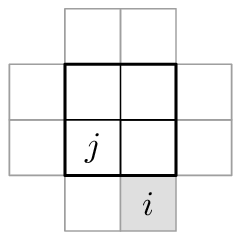}
\caption{\label{fig:i-and-j-2x2-block}A 2$\times$2-block $\Theta_{k}$ showing
all eight positions of a vertex $i\in\partial\Theta_{k}$ on the boundary
of the block in relation to a vertex $j\in\Theta_{k}$ in the block.}
\end{figure}

Thus if $i\in\partial\Theta_{k}$ we let $\Psi_{k}(\sigma,\tau)$
be the coupling of $\pi_{k}(\sigma)$ and $\pi_{k}(\tau)$ that draws
the colouring of $\Theta_{k}$ from the coupling $\Psi$ in Lemma~\ref{lem:main-lemma},
where $Z$ is the boundary colouring obtained from $\sigma_{\partial\Theta_{k}}$
and $Z'$ is obtained from $\tau_{\partial\Theta_{k}}$, and leaves
the colour of the remaining vertices, $V\backslash\Theta_{k}$, unchanged.
That is, if the pair $(\sigma',\tau')$ of colourings are drawn from
$\Psi_{k}(\sigma,\tau)$ then $\sigma'=\sigma$ off $\Theta_{k}$,
$\tau'=\tau$ off $\Theta_{k}$ and the colourings of $\Theta_{k}$
in $\sigma'$ and $\tau'$ are drawn from the coupling $\Psi$ in
Lemma~\ref{lem:main-lemma} (see the proof for details on how to
construct $\Psi$). It is straightforward to verify that this is indeed
a coupling of $\pi_{k}(\sigma)$ and $\pi_{k}(\tau)$. Note that due
to the symmetry of the 2$\times$2-block, with respect to rotation
and mirroring, we can always label the vertices of $\Theta_{k}$ and
$\partial\Theta_{k}$ such that label $z_{1}$ in Figure~\ref{fig:labeling-2x2-block}
represents the discrepancy vertex $i$ on the boundary. Hence we can
make use of Lemma~\ref{lem:main-lemma} to compute upper bounds on
the parameters $\rho_{i,j}^{k}$. We summarise the $\rho_{i,j}^{k}$-values
in the following Corollary of Lemma~\ref{lem:main-lemma}. Note that
due to the symmetry of the block we can assume that vertex $j\in\Theta_{k}$
in the corollary is located in the bottom left corner, as Figure~\ref{fig:i-and-j-2x2-block}
shows.

\begin{corollary}
\label{cor:bounds}Let $\Theta_{k}$ be any 2$\times$2-block, let
$j\in\Theta_{k}$ be any vertex in the block and let $i\in\partial\Theta_{k}$
be a vertex on the boundary of the block. Then
\begin{equation}
\rho_{i,j}^{k}=\max_{(\sigma,\tau)\in S_{i}}p_{j}(\Psi_{k}(\sigma,\tau))<
\begin{cases}
\pone, & \textup{if }i\textup{ and }j\textup{ as in Figure \ref{fig:i-and-j-2x2-block}(a) or (b),}\\
\ptwo, & \textup{if }i\textup{ and }j\textup{ as in Figure \ref{fig:i-and-j-2x2-block}(c) or (h),}\\
\pthree, & \textup{if }i\textup{ and }j\textup{ as in Figure \ref{fig:i-and-j-2x2-block}(e) or (f),}\\
\pfour, & \textup{if }i\textup{ and }j\textup{ as in Figure \ref{fig:i-and-j-2x2-block}(d) or (g).}
\end{cases}
\end{equation}
If $i\notin\partial\Theta_{k}$ is not on the boundary of the block
then $\rho_{i,j}^{k}=0$.
\end{corollary}

We can then use Corollary~\ref{cor:bounds} to prove Theorem~\ref{thm:gridmix}.
The proof of Theorem~ \ref{thm:gridmix} is given here:

\begin{proof}
[Proof of Theorem~\ref{thm:gridmix}]
Let $\alpha_{k,j}=\sum_{i}\rho_{i,j}^{k}$
be the influence on $j$ under $\Theta_{k}$. We need $\alpha_{k,j}$
to be upper bounded by one for each block $\Theta_{k}$ and vertex
$j\in\Theta_{k}$ in order to ensure that $\alpha=\max_{k}\max_{j\in\Theta_{k}}\alpha_{k,j}$
is less than one. Fix any block $\Theta_{k}$ and any vertex $j\in\Theta_{k}$.
A vertex $i\in\partial\Theta_{k}$ on the boundary of the block can
occupy eight different positions on the boundary in relation to $j$
as showed in Figure~\ref{fig:i-and-j-2x2-block}(a)--(h).
Recall that we are working on the torus, and hence every vertex on the boundary of the block will belong to~$G$.
Thus, using
the bounds from Corollary~\ref{cor:bounds} we have
\begin{equation}
\label{eq:alpha-k-j}
\alpha_{k,j}=\sum_{i}\rho_{i,j}^{k}<2(\pone+\ptwo+\pthree+\pfour)=\alphasum.
\end{equation} 
Then $\alpha=\max_{k}\max_{j\in\Theta_{k}}\alpha_{k,j}<\max_{k}\alphasum=\alphasum<1$
and we obtain the stated bound on the mixing time of $\gridscan$
by Theorem~\ref{thm:main-d}.
\end{proof}

We make the following remark.
In the proof of Theorem \ref{thm:gridmix} above, we assume that $G$ is a finite rectangular grid
with toroidal boundary conditions. Hence, every block is a 2$\times$2-sub-grid and each vertex on the
block boundary belongs to $V$. We note that if $G$ is a finite rectangular grid without
toroidal boundary conditions
then some vertices on the boundary $\partial \Theta_k$ of a block $\Theta_k$
might fall outside $G$.
The sum in Equation~(\ref{eq:alpha-k-j}) is over boundary vertices $i$ that do belong to $V$, and hence
the number of terms in this sum is reduced if some boundary vertices do not belong to $V$,
making $\alpha$ smaller.
Furthermore, if $G$ is a non-rectangular region of the grid then a block next to the boundary
might be smaller than 2$\times$2 vertices.
Suppose $\Theta_k$ is a block that is smaller than 2$\times$2 vertices. Then the vertices that
are missing in order to make $\Theta_k$ a full 2$\times$2-block are boundary vertices.
Suppose $i \in \partial \Theta_k$ belongs to $V$ and $i' \in \partial \Theta_k$ does not belong to $V$.
When constructing couplings $\Psi_{k}(\sigma,\tau)$, where $(\sigma,\tau)\in S_i$, we must consider the
vertex $i'$ as ``colourless'', which would decrease the value of $p_{i,j}^k$.
A more rigorous analysis yields that our mixing result with seven~colours and 2$\times$2-blocks holds
for arbitrary finite regions $G$ of the grid.

Of course we have yet to establish a proof of Lemma \ref{lem:main-lemma},
and the rest of this paper will be concerned with this. Our method
of proof uses some ideas of Goldberg, Jalsenius, Martin and Paterson~\cite{ssm_markus}
in so far as it is computer assisted and we will be focusing on minimising
the probability of assigning different colours to vertex $v_{1}$
in the constructed couplings. We will however be required to construct
a coupling on the 2$\times$2 sub-grid, rather than establishing bounds
on the disagreement probability of a vertex adjacent to the initial
discrepancy and then extending this to a coupling on the whole block
recursively. Our approach is similar to the one Achlioptas et al.~\cite{grid_ach}
take, however we do not have the option of constructing an ``optimal''
coupling using a suitable linear program (even when feasible) since
our probabilities will be maximised over all boundary colourings.
The crucial difference between the approaches is that Achlioptas et
al.~\cite{grid_ach} are using path coupling (see Bubley and Dyer~\cite{pathcoupling})
as a proof technique which requires them to bound the expected Hamming
distance between a pair of colourings drawn from a coupling. This in
turn enables them to, for a given boundary colouring, specify an ``optimal''
coupling which minimises Hamming distance. We are, however, required
to bound the influence of $i$ on $j$ for each boundary colouring
and sum over the maximum of these influences. The reason for this
is the inherit maximisation over boundary colourings in the definition
of $\rho_{i,j}^{k}$ as described above. 

Finally it is worth mentioning that providing bounds on the expected
Hamming distance is similar to showing that the influence \emph{of}
a vertex is small and it is known that this condition implies rapid
mixing of a random update Markov chain, see for example Weitz~\cite{dror_combinatorial}.
In a single-site setting the condition ``the influence \emph{of}
a vertex is small'' also implies rapid mixing of systematic scan (Dyer
et al.~\cite{dobrushin_scan}), however, in a block setting this condition is not
sufficient to give rapid mixing of systematic scan (Pedersen~\cite{thesis}),
which is why we need to bound the influence \emph{on} a vertex.

\section{Constructing the coupling by machine\label{sec:coupling}}

In order to prove Lemma~\ref{lem:main-lemma} we will construct a
coupling $\Psi$ of $\pi_{Z}$ and $\pi_{Z'}$ for all pairs of boundary
colourings $Z$ and $Z'$ that are identical on all boundary vertices
but vertex $z_{1}$, on which $Z$ and $Z'$ differ. For each coupling
constructed we verify that the probabilities $p_{v_{i}}(\Psi)$, $i=1,\dots,4$,
are within the bounds of the lemma. The method is well suited to be
carried out with the help of a computer and we have implemented a
program in C to do so. Before stating the proof of Lemma~\ref{lem:main-lemma}
we will discuss how a coupling can be represented by an edge-weighted
complete bipartite graph. We make use of this representation of $\Psi$
in the proof of the lemma.

\subsection{Representing a coupling as a bipartite graph}

Let $S$ be a set of objects and let $W$ be a set of $|S|$ pairs
$(s,w_{s})$ such that $s\in S$ and $w_{s}\geq0$ is a non-negative
value representing the weight of $s$. Each element $s\in S$ is contained
in exactly one of the pairs in $W$. If the value $w_{s}$ is an integer
(which it is in our case) it can be regarded as the multiplicity of
$s$ in a multiset. The set $W$ is referred to as a \emph{weighted
set of $S$}. Let $\pi_{S,W}$ be the distribution on $S$ such that
the probability of $s$ is proportional to $w_{s}$, where $(s,w_{s})$
is a pair in $W$. More precisely, the probability of $s$ in $\pi_{S,W}$
is $\Prob_{\pi_{S,W}}(s)=w_{s}/\sum_{(t,w_{t})\in W}w_{t}$. For example,
let $W$ be a weighted set of $S$ and let $S'\subseteq S$ be a subset
of $S$. Assume the weight $w_{s}=0$ if $s\in S\backslash S'$ and
$w_{s}=k$ if $s\in S'$, where $k>0$ is a positive constant. Then
$\pi_{S,W}$ is the uniform distribution on $S'$.

The reason for introducing the notion of a weighted set is that it
can be used when specifying a coupling of two distributions. Let $S$
be a set and let $W$ and $W'$ be two weighted sets of $S$ such
that the sum of the weights in $W$ equals the sum of the weights
in $W'$. Let $w_{\text{tot}}$ denote this sum. That is, $w_{\text{tot}}=\sum_{(s,w_{s})\in W}w_{s}=\sum_{(s',w_{s'}')\in W'}w_{s'}'$.
The two weighted sets $W$ and $W'$ define two distributions $\pi_{S,W}$
and $\pi_{S,W'}$ on $S$. We want to specify a coupling $\Psi$ of
$\pi_{S,W}$ and $\pi_{S,W'}$. Let $K_{|S|,|S|}$ be an edge-weighted
complete bipartite graph with vertex sets $W$ and $W'$. That is,
for each pair $(s,w_{s})\in W$ there is an edge to every pair in
$W'$. Every edge $e$ of $K_{|S|,|S|}$ has a weight $w_{e}\geq0$
such that the following condition holds. Let $(s,w_{s})$ be any pair
in $W\cup W'$ and let $E$ be the set of all $|S|$ edges incident
to $(s,w_{s})$. Then $\sum_{e\in E}w_{e}=w_{s}$. It follows that
the sum of the edge weights of all $|S|^{2}$ edges in $K_{|S|,|S|}$
equals $w_{\text{tot}}$, the sum of the weights in $W$ (and $W'$).
The idea is that $K_{|S|,|S|}$ represents a coupling $\Psi$ of $\pi_{S,W}$
and $\pi_{S,W'}$. In order to draw a pair of elements from $\Psi$
we randomly select an edge $e$ in $K_{|S|,|S|}$ proportional to
its weight. The endpoints of $e$ represent the elements in $S$ drawn
from $\pi_{S,W}$ and $\pi_{S,W'}$. More precisely, the probability
of choosing edge $e$ in $K_{|S|,|S|}$ with weight $w_{e}$ is $w_{e}/w_{\text{tot}}$.
If edge $e=((s,w_{s}),(s',w_{s'}'))$ is chosen it means that we have
drawn $s$ from $\pi_{S,W}$ and $s'$ from $\pi_{S,W'}$, the marginal
distributions of $\Psi$.

The bipartite graph representation of a coupling will be used when
we construct couplings of colourings of 2$\times$2-blocks in the
proof of Lemma~\ref{lem:main-lemma}.

\subsection{The proof of Lemma~\ref{lem:main-lemma}}

Here is the proof of Lemma~\ref{lem:main-lemma}:
\begin{proof}
[Proof of Lemma~\ref{lem:main-lemma}]Fix two colourings $Z$ and
$Z'$ of the boundary that differ on vertex $z_{1}$. Let $c$ be
the colour of vertex $z_{1}$ in $Z$ and let $c'\neq c$ be the colour
of $z_{1}$ in $Z'$. Let $C_{Z}$ and $C{}_{Z'}$ be the two sets
of proper 7-colourings of the block that agree with $Z$ and $Z'$,
respectively. Let $C^{+}$ be the set of all 7-colourings of the block.
Let $W_{Z}$ and $W_{Z'}$ be two weighted sets of $C^{+}$. The weights
are assigned as follows.
\begin{itemize}
\item For the pair $(\sigma,w_{\sigma})\in W_{Z}$ let the weight $w_{\sigma}=|C_{Z'}|$
if $\sigma\in C_{Z}$, otherwise let $w_{\sigma}=0$.
\item For the pair $(\sigma,w_{\sigma})\in W_{Z'}$ let the weight $w_{\sigma}=|C_{Z}|$
if $\sigma\in C_{Z'}$, otherwise let $w_{\sigma}=0$. 
\end{itemize}
It follows from the assignment of the weights that the distribution
$\pi_{C^{+},W_{Z}}$ is the uniform distribution on $C_{Z}$. That
is, $\pi_{C^{+},W_{Z}}=\pi_{Z}$. Similarly, $\pi_{C^{+},W_{Z'}}$
is the uniform distribution $\pi_{Z'}$ on $C_{Z'}$. Note that the
sum of the weights is $|C_{Z}||C_{Z'}|$ in both $W_{Z}$ and $W_{Z'}$.
Then a coupling $\Psi$ of $\pi_{C^{+},W_{Z}}$ and $\pi_{C^{+},W_{Z'}}$
can be specified with an edge-weighted complete bipartite graph $K=K_{|C^{+}|,|C^{+}|}$.
For a given valid assignment of the weights of the edges of $K$,
making $K$ represent a coupling $\Psi$, we can compute the probabilities
of having a mismatch on a vertex $v_{i}$ of the block when two colourings
are drawn from $\Psi$. Let $E$ be the set of all edges $e=((\sigma,w_{\sigma}),(\sigma',w_{\sigma'}'))$
in $K$ such that $\sigma$ and $\sigma'$ differ on vertex $v_{i}$.
Then $p_{v_{i}}(\Psi)=\sum_{e\in E}w_{e}/|C_{Z}||C_{Z'}|$.

In order to obtain sufficiently small upper bounds on $p_{v_{i}}(\Psi)$
for the four vertices $v_{1},\dots,v_{4}$ in the block we would like
to assign weights to the edges of $K$ such that much weight is assigned
to edges between colourings that agree on many vertices in the block.
In general it is not clear exactly how to assign weights to the edges.
For instance, if we assign too much weight to edges between colourings
that are identical on vertex $v_{2}$ we might not be able to assign
as much weight as we would like to on edges between colourings that
are identical on vertex $v_{4}$. Thus, the probability of having
a mismatch on $v_{4}$ would increase. Intuitively a good strategy
would be to assign as much weight as possible to edges between colourings
that are identical on the whole block. This implies that we try to
assign as much weight as possible to edges between colourings that
are identical on vertex $v_{1}$, the vertex adjacent to the discrepancy
vertex $z_{1}$ on the boundary. If there is a mismatch on vertex
$v_{1}$ it should be a good idea to assign as much weight as possible
to edges between colourings that are identical on the whole block
apart from vertex $v_{1}$. This idea leads to a heuristic in which
the assignment of the edge weights is divided into three phases. The
exact procedure is described as follows.

In phase one we match identical colourings. For all colourings $\sigma\in C^{+}$
of the block the edge $e=((\sigma,w_{\sigma}),(\sigma,w_{\sigma}'))$
in $K$ will be given weight $w_{e}=\min(w_{\sigma},w_{\sigma}')$.
That is, we maximise the probability of drawing the same colouring
$\sigma$ from both $\pi_{C^{+},W_{Z}}$ and $\pi_{C^{+},W_{Z'}}$.

For the following two phases we define an ordering of the colourings
in $C^{+}$. We order the colourings lexicographically with respect
to the vertex order $v_{3}$, $v_{2}$, $v_{4}$, $v_{1}$. That is,
if the seven colours are $1,\dots,7$ the colouring of $v_{3}$, $v_{2}$,
$v_{4}$, $v_{1}$ will start with 1, 1, 1, 1, respectively. The next
colouring will be 1, 1, 1, 2, and so on. This ordering of colourings
in $C^{+}$ carries over to an ordering of the pairs in $W_{Z}$ and
$W_{Z'}$. That is, we order the pairs $(\sigma,w_{\sigma})$ in $W_{Z}$
with respect to the lexicographical ordering of $\sigma$. Similarly
we order the pairs in $W_{Z'}$. This ordering of the pairs will be
important in the next two phases. It provides some control of how
colourings are being paired up in terms of the assignment of the weights
on edges between pairs. Edges will be considered with respect to this
ordering because choosing an arbitrary ordering of the edges would
not necessarily result in probabilities $p_{v_{i}}(\Psi)$ that would
be within the bounds of the lemma.

In the second phase we ignore the colour of vertex $v_{1}$ and match
colourings that are identical on all of the remaining three vertices
$v_{2}$, $v_{3}$ and $v_{4}$. More precisely, for each pair $(\sigma,w_{\sigma})\in W_{Z}$,
considered in the ordering explained above, we consider the edges
$e=((\sigma,w_{\sigma}),(\sigma',w_{\sigma'}'))$ where $\sigma$
and $\sigma'$ are identical on all vertices but $v_{1}$. The edges
are considered in the ordering of the second component $(\sigma',w_{\sigma'}')\in W_{Z'}$.
We assign as much weight as possible to $e$ such that the total weight
on edges incident to $(\sigma,w_{\sigma})\in W_{Z}$ does not exceed
$w_{\sigma}$ and such that the total weight on edges incident to
$(\sigma',w_{\sigma'}')\in W_{Z'}$ does not exceed $w_{\sigma'}'$.
Note that in the lexicographical ordering of the colourings, vertex
$v_{1}$ is the least significant vertex and therefore the ordering
provides some level of control of pairing up colourings that are similar
on the remaining three vertices. It turns out that the resulting coupling
is sufficiently good for proving the lemma.

In the third and last phase we assign the remaining weights on the
edges. As in phase two, for each pair $(\sigma,w_{\sigma})\in W_{Z}$
we consider the edges $e=((\sigma,w_{\sigma}),(\sigma',w_{\sigma'}'))$.
The pairs and edges are considered in accordance with the ordering
explained above. The difference between the second and third phase
is that now we do not have any restrictions on the colourings $\sigma$
and $\sigma'$. We assign as much weight as possible to $e$ such
that the total weight on edges incident to $(\sigma,w_{\sigma})\in W_{Z}$
does not exceed $w_{\sigma}$ and such that the total weight on edges
incident to $(\sigma',w_{\sigma'}')\in W_{Z'}$ does not exceed $w_{\sigma'}'$.
After phase three we have assigned all weights to the edges of $K$
and hence $K$ represents a coupling $\Psi$ of $\pi_{Z}$ and $\pi_{Z'}$.

From $K$ we compute the probabilities $p_{v_{1}}(\Psi)$, $p_{v_{2}}(\Psi)$,
$p_{v_{3}}(\Psi)$ and $p_{v_{4}}(\Psi)$ as described above. We have
written a C-program which loops through all colourings $Z$ and $Z'$
of the boundary of the block and constructs the bipartite graph $K$
as described above. For each boundary the probabilities $p_{v_{1}}(\Psi)$,
$p_{v_{2}}(\Psi)$, $p_{v_{3}}(\Psi)$ and $p_{v_{4}}(\Psi)$ are
successfully verified to be within the bounds of the lemma. For details
on the C-program, see \scanURL.
\end{proof}

\section{Partial results for $6$-colourings of the grid \label{sec:Partial-results}}

In previous sections we have seen that systematic scan on the grid
using 2$\times$2-blocks and seven colours mixes rapidly. An immediate
question is whether we can do better and show rapid mixing with six
colours. This matter will be discussed in this section and we will
show that, even with bigger block sizes (up to 3$\times$3), it is
not possible to show rapid mixing using the technique of this paper.
More precisely, we will establish lower bounds on the parameter $\alpha$
for 2$\times$2-blocks, 2$\times$3-blocks and 3$\times$3-blocks.
All three lower bounds are greater than one and hence we cannot make
use of Theorem~\ref{thm:main-d} to show rapid mixing.

\subsection{Establishing lower bounds for 2$\times$2 blocks}

We start by examining the 2$\times$2-block again but this time with
six colours. Lemma~\ref{lem:main-lemma} provides upper bounds (under
any colourings of the boundary) on the probabilities of having discrepancies
at each of the four vertices of the block when two 7-colourings are
drawn from the specified coupling. For six colours we will show lower
bounds on these probabilities under any coupling and a specified pair
of boundary colourings. Once again, let $v_{1},\dots,v_{4}$ be the
four vertices in a 2$\times$2-block and let $z_{1},\dots,z_{8}$
be the boundary vertices of the block and let the labeling be as in
Figure~\ref{fig:labeling-2x2-block}. Let $Z$ and $Z'$ be any two
$6$-colourings of the boundary vertices that assign the same colour
to each vertex except for $z_{1}$. Let $\pi_{Z}$ and $\pi_{Z'}$
be the uniform distributions on the sets of proper $6$-colourings
of the block that agree with $Z$ and $Z'$, respectively. Let $\Psi_{v_{k}}^{\textup{min}}(Z,Z')$
be a coupling of $\pi_{Z}$ and $\pi_{Z'}$ that minimises $p_{v_{k}}(\Psi)$.
That is, $p_{v_{k}}(\Psi)\geq p_{v_{k}}(\Psi_{v_{k}}^{\textup{min}}(Z,Z'))$
for all couplings $\Psi$ of $\pi_{Z}$ and $\pi_{Z'}$. Also let
$p_{v_{k}}^{\textup{low}}=\max_{Z,Z'}p_{v_{k}}(\Psi_{v_{k}}^{\textup{min}}(Z,Z'))$.
We can hence say that there exist two 6-colourings $Z$ and $Z'$
of the boundary of a 2$\times$2 block, that assign the same colour
to each vertex except for $z_{1}$, such that $p_{v_{k}}(\Psi)\geq p_{v_{i}}^{\textup{low}}$
for any coupling $\Psi$ of $\pi_{Z}$ and $\pi_{Z'}$. We have the
following lemma, which is proved by computation.

\begin{lemma}
\label{lem:6-colours-probs-2x2}Consider 6-colourings of the 2$\times$2-block
in Figure~\ref{fig:labeling-2x2-block}. Then $p_{v_{1}}^{\textup{low}}\geq\ponesmall$,
$p_{v_{2}}^{\textup{low}}\geq\ptwosmall$, $p_{v_{3}}^{\textup{low}}\geq\pthreesmall$
and $p_{v_{4}}^{\textup{low}}\geq\pfoursmall$.
\end{lemma}
\begin{proof}
Fix one vertex $v_{k}$ in the block and fix two colourings $Z$ and
$Z'$ of the boundary of the block that differ only on the colour
of vertex $z_{1}$. Let $C_{Z}$ and $C_{Z'}$ be the two sets of
proper 6-colourings of the block that agree with $Z$ and $Z'$, respectively.
For $c=1,\dots,6$ let $n_{c}$ be the number of colourings in $C_{Z}$
in which vertex $v_{k}$ is assigned colour $c$. Similarly let $n_{c}'$
be the number of colourings in $C_{Z'}$ in which vertex $v_{k}$
is assigned colour $c$. It is clear that the probability that $v_{k}$
is assigned colour $c$ in a colouring $\sigma'$ drawn from $\pi_{Z}$
is $\Prob_{\pi_{Z}}(\sigma'_{v_{k}}=c)=n_{c}/|C_{Z}|$. For $c=1,\dots,6$
define $m_{c}=n_{c}|C_{Z'}|$, $m_{c}'=n_{c}'|C_{Z}|$ and $M=|C_{Z}||C_{Z'}|$.
It follows that $\Prob_{\pi_{Z}}(\sigma'_{v_{k}}=c)=m_{c}/M$ and
$\Prob_{\pi_{Z'}}(\tau'_{v_{k}}=c)=m_{c}'/M$, where $\sigma'$ and
$\tau'$ are colourings drawn from $\pi_{Z}$ and $\pi_{Z'}$, respectively.
Observe that the quantities $m_{c}$, $m_{c}'$ and $M$ can be easily
computed for a given pair of boundary colourings.

Now let $\Psi$ be any coupling of $\pi_{Z}$ and $\pi_{Z'}$. It
is easy to see that the probability that vertex $v_{k}$ is coloured
$c$ in both colourings drawn from $\Psi$ can be at most $\min(m_{c},m_{c}')/M$.
Therefore, the probability of drawing two colourings from $\Psi$
such that the colour of vertex $v_{k}$ is the same in both colourings
is at most $\sum_{c=1,\dots,6}\min(m_{c},m_{c}')/M$, and the probability
of assigning different colours to vertex $v_{k}$ is at least $p_{v_{k}}(\Psi)\geq1-\sum_{c=1,\dots,6}\min(m_{c},m_{c}')/M$.
We have successfully verified the bounds in the statement of the lemma
by maximising the lower bound on $p_{v_{k}}(\Psi)$ over all boundary
colourings $Z$ and $Z'$ for each vertex $v_{k}$ in the block. The
computations are carried out with the help of a computer program written
in C. For details on the program, see \scanURL.
\end{proof}
For seven colours, Corollary~\ref{cor:bounds} makes use of Lemma~\ref{lem:main-lemma}
to establish upper bounds on the influence parameters $\rho_{i,j}^{k}$.
These parameters are used in the proof of Theorem~\ref{thm:gridmix}
to obtain an upper bound on the parameter $\alpha$. The upper bound
on $\alpha$ is shown to be less than one which implies rapid mixing
for seven colours when applying Theorem~\ref{thm:main-d}. We can
use Lemma~\ref{lem:6-colours-probs-2x2} to obtain lower bounds on
the influence parameters $\rho_{i,j}^{k}$ by completing the coupling
in a way analogous to the coupling in Corollary~\ref{cor:bounds}.
This in turn will result in a lower bound on the parameter $\alpha$
that is greater than one. That is, following the proof of Theorem~\ref{thm:gridmix}
and making use of Lemma~\ref{lem:6-colours-probs-2x2}, a lower bound
on $\alpha$ will be
\begin{equation}
\alpha\geq2(\ponesmall+\ptwosmall+\pthreesmall+\pfoursmall)=\alphasumsmall>1.
\end{equation}
Hence we fail to show rapid mixing of systematic scan with six colours
using 2$\times$2-blocks.

\subsection{Bigger blocks}

\begin{figure}
\centering
(a)\includegraphics[scale=0.8]{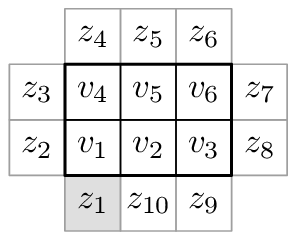}
\hspace{3mm}
(b)\includegraphics[scale=0.8]{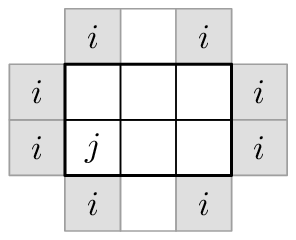}
\hspace{3mm}
(c)\includegraphics[scale=0.8]{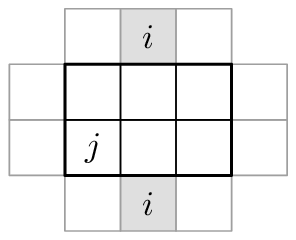}
\caption{\label{fig:2x3-block}(a) General labeling of the vertices in a 2$\times$3-block
$\Theta_{k}$ and the vertices $\partial\Theta_{k}$ on the boundary
of the block. (b)--(c) All ten positions of a vertex $i\in\partial\Theta_{k}$
on the boundary of the block in relation to a vertex $j\in\Theta_{k}$
in the corner of the block.}
\end{figure}

We failed to show rapid mixing of systematic scan with six colours
and 2$\times$2-blocks and we will now show that increasing the block
size to both 2$\times$3 and 3$\times$3 will not be sufficient either.
Lemma~\ref{lem:6-colours-probs-2x3} below considers 2$\times$3-blocks
and is analogous to Lemma~\ref{lem:6-colours-probs-2x2}. We make
use of the same notation as for Lemma~\ref{lem:6-colours-probs-2x2},
only the block is bigger and the labeling of the vertices is different
(see Figure~\ref{fig:2x3-block}(a)).
Lemma~\ref{lem:6-colours-probs-2x3} is proved by computation in
the same way as Lemma~\ref{lem:6-colours-probs-2x2}. For details
on the C-program used in the proof, see \scanURL.

\begin{lemma}
\label{lem:6-colours-probs-2x3}Consider 6-colourings of the 2$\times$3-block
in Figure~\ref{fig:2x3-block}(a). Then $p_{v_{1}}^{\textup{low}}\geq\ponemedium$,
$p_{v_{3}}^{\textup{low}}\geq\pthreemedium$, $p_{v_{4}}^{\textup{low}}\geq\pfourmedium$
and $p_{v_{6}}^{\textup{low}}\geq\psixmedium$.
\end{lemma}

We will now use Lemma~\ref{lem:6-colours-probs-2x3} to show that
$\alpha>1$ for 2$\times$3 blocks. Let $\Theta_{k}$ be any 2$\times$3-block
and let $j\in\Theta_{k}$ be a vertex in a corner of the block. A
vertex $i\in\partial\Theta_{k}$ on the boundary of the block can
occupy ten different positions on the boundary in relation to $j$.
See Figure~\ref{fig:2x3-block}(b) and~(c). We can again determine
lower bounds on the influences $\rho_{i,j}^{k}$ of $i$ on $j$ under
$\Theta_{k}$ from Lemma~\ref{lem:6-colours-probs-2x3}. However,
Lemma~\ref{lem:6-colours-probs-2x3} provides lower bounds on $\rho_{i,j}^{k}$
only when $i\in\partial\Theta_{k}$ is adjacent to a corner vertex
of the block, as in Figure~\ref{fig:2x3-block}(b). If $i$ is located
as in Figure~\ref{fig:2x3-block}(c) we do not know more than that
$\rho_{i,j}^{k}$ is bounded from below by zero. Nevertheless, the
lower bound on $\alpha$ exceeds one. Let $\alpha_{k,j}=\sum_{i}\rho_{i,j}^{k}$
be the influence on $j$ under $\Theta_{k}$. Following the proof
of Theorem~\ref{thm:gridmix} and using the lower bounds in Lemma~\ref{lem:6-colours-probs-2x3}
we have
\begin{eqnarray}
\alpha_{k,j} & = & \sum_{i\textup{ in Fig. \ref{fig:2x3-block}(b)}}\rho_{i,j}^{k}+\sum_{i\textup{ in Fig. \ref{fig:2x3-block}(c)}}\rho_{i,j}^{k}\notag\\
& \geq & 2(\ponemedium+\pthreemedium+\pfourmedium+\psixmedium)=\alphasummedium,
\end{eqnarray}
where we set the lower bound on the second sum to zero. Now,
\begin{equation}
\alpha=\max_{k}\max_{j\in\Theta_{k}}\alpha_{k,j}\geq\alphasummedium>1.
\end{equation}
Hence we cannot use Theorem~\ref{thm:main-d} to show rapid mixing
of systematic scan with six colours and 2$\times$3-blocks. It is
interesting to note that considering 2$\times$3-blocks was sufficient
for Achlioptas et~al.~\cite{grid_ach} to prove mixing of a random
update Markov chain for sampling 6-colourings of the grid. 

Lastly, we increase the block size to 3$\times$3 and show that a
lower bound on $\alpha$ is still greater than one. We have the following
lemma which is proved by computation in the same way as Lemmas~\ref{lem:6-colours-probs-2x2}
and~\ref{lem:6-colours-probs-2x3}. For details on the C-program
used in the proof see \scanURL.

\begin{figure}
\centering
(a)\includegraphics[scale=0.8]{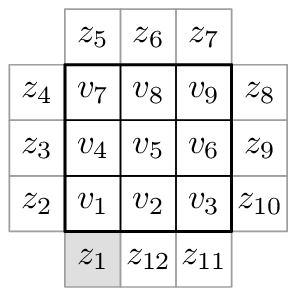}
\hspace{2mm}
(b)\includegraphics[scale=0.8]{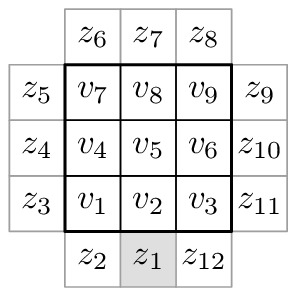}
\hspace{2mm}
(c)\includegraphics[scale=0.8]{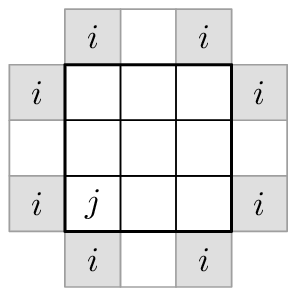}
\hspace{2mm}
(d)\includegraphics[scale=0.8]{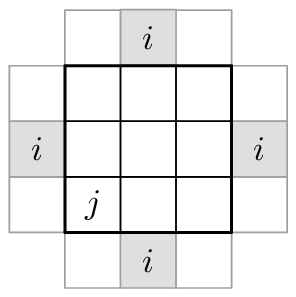}
\caption{\label{fig:3x3-block}(a)--(b) General labeling of the vertices in
a 3$\times$3-block $\Theta_{k}$ and two different labellings of
the vertices $\partial\Theta_{k}$ on the boundary of the block. The
discrepancy vertex on the boundary has label $z_{1}$. (b)--(c) All
twelve positions of a vertex $i\in\partial\Theta_{k}$ on the boundary
of the block in relation to a vertex $j\in\Theta_{k}$ in the corner
of the block.}
\end{figure}

\begin{lemma}
\label{lem:6-colours-probs-3x3}For 6-colourings of the 3$\times$3-block
with vertices labeled as in Figure~\ref{fig:3x3-block}(a) we have
$p_{v_{1}}^{\textup{low}}\geq\ponecornerbig$, $p_{v_{3}}^{\textup{low}}\geq\pthreecornerbig$,
$p_{v_{7}}^{\textup{low}}\geq\psevencornerbig$ and $p_{v_{9}}^{\textup{low}}\geq\pninecornerbig$.
Furthermore, for 6-colourings of the 3$\times$3-block in Figure~\ref{fig:3x3-block}(b)
we have $p_{v_{1}}^{\textup{low}}\geq\ponemiddlebig$, $p_{v_{3}}^{\textup{low}}\geq\pthreemiddlebig$,
$p_{v_{7}}^{\textup{low}}\geq\psevenmiddlebig$ and $p_{v_{9}}^{\textup{low}}\geq\pninemiddlebig$.
\end{lemma}

Note that Lemma~\ref{lem:6-colours-probs-3x3} provides lower bounds
on the probabilities of having a mismatch on a corner vertex of the
block when the discrepancy vertex on the boundary (labeled $z_{1})$
is adjacent to a corner vertex (Figure~\ref{fig:3x3-block}(a)) and
adjacent to a middle vertex (Figure~\ref{fig:3x3-block}(b)). Let
$\Theta_{k}$ be any 3$\times$3-block and let $j\in\Theta_{k}$ be
a vertex in a corner of the block. A vertex $i\in\partial\Theta_{k}$
on the boundary of the block can occupy twelve different positions
on the boundary in relation to $j$. See Figure~\ref{fig:3x3-block}(c)
and~(d). Analogous to Corollary~\ref{cor:bounds} lower bounds on
the influences $\rho_{i,j}^{k}$ of $i$ on $j$ under $\Theta_{k}$
can be determined from Lemma~\ref{lem:6-colours-probs-3x3}. Let
$\alpha_{k,j}=\sum_{i}\rho_{i,j}^{k}$ be the influence on $j$ under
$\Theta_{k}$. Following the proof of Theorem~\ref{thm:gridmix}
and using the lower bounds in Lemma~\ref{lem:6-colours-probs-3x3}
we have
\begin{eqnarray}
\alpha_{k,j} & = & \sum_{i\textup{ in Fig. \ref{fig:3x3-block}(c)}}\rho_{i,j}^{k}+\sum_{i\textup{ in Fig. \ref{fig:3x3-block}(d)}}\rho_{i,j}^{k}\notag\\
 & \geq & 2(\ponecornerbig+\pthreecornerbig+\psevencornerbig+\pninecornerbig)+\notag\\
 & & (\ponemiddlebig+\pthreemiddlebig+\psevenmiddlebig+\pninemiddlebig)=\alphasumbig.
\end{eqnarray}
Thus, $\alpha=\max_{k}\max_{j\in\Theta_{k}}\alpha_{k,j}\geq\alphasumbig>1$.
Hence, we cannot use Theorem~\ref{thm:main-d} to show rapid mixing
of systematic scan with six colours and 3$\times$3-blocks. 

A natural question is whether we can show rapid mixing using even
bigger blocks. It seems possible to do this although the computations
rapidly become intractable as the block size increases. Already with
a 3$\times$3-block the number of boundary colourings we need to consider
(after removing isomorphisms) is in excess of $10^{6}$ and for each
boundary colouring there are more than $10^{7}$ colourings of the
block to consider. In addition to simply generating the distributions
on colourings of the block, the time it would take to actually construct
the required couplings, as we did in the proof of Lemma~\ref{lem:main-lemma},
would also increase. Finally when using a larger block size, different
positions of vertex $j$ in the block need to be considered whereas
we could make use of to the symmetry of the 2$\times$2-block to only
consider one position of vertex $j$ in the block. If different positions
of $j$ have to be considered this has to be captured in the construction
of the coupling and would likely require more computations. The conclusion
is that in order to show rapid mixing for six colours of systematic
scan on the grid we would most likely have to rely on a different
approach than the one presented in this paper.

\bibliographystyle{plain}

\end{document}